\documentclass[11pt]{article}
\usepackage{amsmath}
\usepackage{amsthm}
\usepackage{graphicx}
\usepackage{amssymb}
\usepackage{psfrag}
\title{Exotic projective structures and \\ 
quasifuchsian spaces II} 
\author{Kentaro Ito}
\theoremstyle{plain}
  \newtheorem{thm}{Theorem}[section]
  \newtheorem{cor}[thm]{Corollary}
  \newtheorem{lem}[thm]{Lemma}
  \newtheorem{prop}[thm]{Proposition}
    
\theoremstyle{definition}
  \newtheorem{defn}[thm]{Definition}
\theoremstyle{remark}
  \newtheorem*{rem}{Remark}

\newcommand{\lam}{\lambda} 
\newcommand{\Lam}{\Lambda } 
 
\newcommand{\C}{\mathbb C} 
\newcommand{\Z}{\mathbb Z} 
\newcommand{\N}{\mathbb N}  
\newcommand{\bd}{\partial} 
 
\newcommand{\hatC}{\widehat{\mathbb C}}   
\newcommand{\what}{\widehat}
\newcommand{\wh}{\widehat}
\newcommand{\wtil}{\widetilde}
\newcommand{\wt}{\widetilde}

\newcommand{\hol}{\mathit{hol}}
\newcommand{\Gr}{{\mathrm{Gr}}}

\newcommand{\QF}{\mathcal{QF}}  
\newcommand{\Q}{\mathcal{Q}}  
\newcommand{\id}{\mathit{id}}
\newcommand{\ov}{\overline}
\newcommand{\psl}{{\mathrm{PSL}}_2(\mathbb C)}
\newcommand{\mln}{{\mathcal{ML}}_{\mathbb N}} 

\newcommand{\inte}{\mathrm{int}\,}
\newcommand{\ul}{\underline}
\newcommand{\LY}{\what{\Lambda}_{Y_\infty}}
\newcommand{\LZ}{\what{\Lambda}_{Z_\infty}}
\newcommand{\LT}{\what{\Lambda}_{T}}
\newcommand{\fix}{\mathrm{fix}\,}
\newcommand{\rn}{{\cal N}}

\newcommand{\fty}{\infty}
\newcommand{\sm}{\setminus}
\newcommand{\hto}{\overset{\mathrm{H}}{\longrightarrow}}

\begin{document}

\maketitle 

\begin{abstract} 
Let $P(S)$ be the space of projective structures on a 
closed surface $S$ of genus $g >1$ 
and let $Q(S)$ be the subset of $P(S)$ 
of projective structures with quasifuchsian holonomy. 
It is known that $Q(S)$ consists of 
infinitely many connected components. 
In this paper, we will show that the closure of 
any exotic component of $Q(S)$ is not 
a topological manifold with boundary 
and that any two components of $Q(S)$ have 
intersecting closures. 
\end{abstract} 

\section{Introduction} 
Let $S$ be an oriented closed surface of genus $g>1$. 
A projective structure on $S$ is a $(G,X)$-structure, 
where $X$ is the Riemann sphere $\hatC$ and 
$G=\psl$ is the group of projective automorphisms of $\hatC$. 
We consider the space of marked projective structures 
$P(S)$ on $S$ and its 
open subset $Q(S)$ of projective 
structures with quasifuchsian holonomy. 
It is known that $Q(S)$ consists 
of infinitely many connected components. 
The aim of this paper is to study how 
components of $Q(S)$ lies in $P(S)$, 
especially how these components bump or self-bump. 
Here we say that components $\Q,\,\Q'$ of $Q(S)$ {\it bump} 
if they have intersecting closures 
and that a component $\Q$ 
{\it self-bumps} if there is a point 
$Y \in \bd \Q$ such that 
for any sufficiently small neighborhood $U$ of $Y$ 
the intersection $U \cap \Q$ is disconnected.

Now let $R(S)$ be the set of conjugacy classes 
of representations $\rho:\pi_1(S) \to \psl$ 
and $\QF \subset R(S)$ the subset
of faithful representations with 
quasifuchsian images. 
It is a result of Hejhal \cite{He} that 
the holonomy map $\hol:P(S) \to R(S)$, 
taking a projective structure to its holonomy representation, 
is a local homeomorphism. 
Therefore, 
studying how $Q(S)=\hol^{-1}(\QF)$ lies in $P(S)$ 
is closely related to studying how the 
quasifuchsian space $\QF$ 
lies in the representation space $R(S)$. 
It is known by Goldman \cite{Go} that 
the set of connected components of $Q(S)$ are classified by 
the set $\mln=\mln(S)$ of integral points of measured laminations; 
see \S 2.3--2.4. 
We denote by $\Q_\lam$ the component of 
$Q(S)$ associated to $\lam \in \mln$, 
where $\Q_0$ is the component of 
standard projective structures. 
Here we say that $Y \in Q(S)$ is {\it standard} 
if its developing map is injective; otherwise it is {\it exotic}. 
One of the most important result on the distribution of components of $Q(S)$ 
in $P(S)$ is obtained by McMullen; 
see Appendix A in \cite{Mc}: 

\begin{thm}[McMullen]
There exists a sequence 
of exotic projective structures 
which converges to a point in the relative boundary 
$\bd \Q_0=\ov{\Q_0} -\Q_0$ 
of the standard component $\Q_0$ in $P(S)$. 
\end{thm}

In fact, given non-zero $\lam \in \mln$, 
McMullen obtained in \cite{Mc} a both-sides infinite convergent sequence 
$$
Y_n \to Y_\fty \in \bd \Q_0 \quad (|n| \to \fty)
$$ 
associated to $\lam$ 
by using the method developed by Anderson and Canary \cite{AnCa}.  
In addition, we showed in \cite{It1} that the projective structures 
$Y_n$ above are contained in the exotic component $\Q_\lam$ for all large $|n|$, 
and thus obtained the following: 

\begin{thm}[Theorem A in \cite{It1}]
For any non-zero $\lam \in \mln$, we have 
$\ov{\Q_0} \cap \ov{\Q_\lam} \ne \emptyset$. 
\end{thm}

Now let $\mu$ be a non-zero element of 
$\mln$ which has no parallel component in common with $\lam$. 
Then one can obtain the grafting 
$Z_\fty=\Gr_\mu(Y_\fty)$ of $Y_\fty$ along $\mu$ (see \S 2.4),  
which lies in the boundary of $\Q_\mu$. 
Since the the map $\hol$ is a local homeomorphism, 
there is a both-sides infinite convergent sequence 
$$
Z_n \to Z_\fty \in \bd \Q_\mu \quad  (|n| \to \fty)
$$
which satisfies 
$\hol(Z_n)=\hol(Y_n)$ for all large $|n|$. 
Although the sequences $\{Y_n\}_{n \gg 0}$ 
and $\{Y_n\}_{n \ll 0}$ are contained in the same component $\Q_\lam$, 
the sequences  $\{Z_n\}_{n \gg 0}$ 
and $\{Z_n\}_{n \ll 0}$ are not necessarily contained in the same component. 
In fact, the main theorem in this paper (Theorem 1.3 below) 
states that these sequences are contained in components 
corresponding to elements 
$(\lam,\mu)_\sharp, \,(\lam, \mu)_\flat$ in $\mln$, respectively, which are  
defined in \S 2.5. See Figure 1. 
We just remark here that 
$(\lam,\mu)_\sharp \ne (\lam, \mu)_\flat$ if and only if 
$i(\lam,\mu) \ne 0$, where $i(\lam,\mu)$ is 
the geometric intersection number of $\lam$ and $\mu$. 

\begin{thm}
In the same notation as above, we have 
$\{Z_n\}_{n \gg 0} \subset \Q_{(\lambda, \mu)_{\sharp}}$ and  
$\{Z_n\}_{n \ll 0} \subset \Q_{(\lambda, \mu)_{\flat}}$.  
\end{thm} 
\bigskip
\begin{figure}[htbp]
\begin{center}
\includegraphics{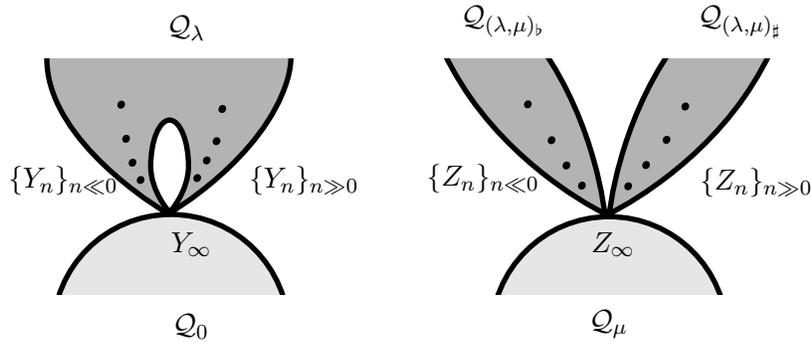}
\end{center}
\vspace{-5cm}
\begin{center}
\unitlength 0.1in
\begin{picture}( 49.5000, 19.0000)( 17.1000,-25.6000)
\put(29.3000,-9.5000){\makebox(0,0)[lb]{$\Q_\lam$}}%
\put(44.7000,-9.3000){\makebox(0,0)[lb]{$\Q_{(\lam,\mu)_\flat}$}}%
\put(29.6000,-25.1000){\makebox(0,0)[lb]{$\Q_0$}}%
\put(51.4000,-25.1000){\makebox(0,0)[lb]{$\Q_\mu$}}%
\put(57.1000,-9.4000){\makebox(0,0)[lb]{$\Q_{(\lam,\mu)_\sharp}$}}%
\put(21.0000,-17.6000){\makebox(0,0)[lb]{$\{Y_n\}_{n \ll 0}$}}%
\put(33.6000,-17.6000){\makebox(0,0)[lb]{$\{Y_n\}_{n \gg 0}$}}%
\put(42.9000,-17.5000){\makebox(0,0)[lb]{$\{Z_n\}_{n \ll 0}$}}%
\put(57.2000,-17.7000){\makebox(0,0)[lb]{$\{Z_n\}_{n \gg 0}$}}%
\put(29.5000,-20.8000){\makebox(0,0)[lb]{$Y_\fty$}}%
\put(51.5000,-20.8000){\makebox(0,0)[lb]{$Z_\fty$}}%
%
\special{pn 8}%
\special{pa 1710 660}%
\special{pa 6660 660}%
\special{pa 6660 2560}%
\special{pa 1710 2560}%
\special{pa 1710 660}%
\special{ip}%
\end{picture}%
\end{center}
\caption{Sequences $\{Y_n\}_{|n| \gg 0}$ and 
$\{Z_n\}_{|n| \gg 0}$.}
\end{figure}

As consequences of Theorem 1.3, we obtain the following results 
on the distribution of components of $Q(S)$; see Theorems 4.1, 4.3 and 4.4. 

\begin{cor}
\begin{enumerate}  
\item  
For any non-zero $\lambda \in \mln$,
$\Q_\lam$ self-bumps. 
\item 
For any $\lambda, \mu \in \mln$, we have
$\overline{\Q_{\lambda}} \cap \overline{\Q_{\mu}} \ne\emptyset$.
\item 
For any non-zero $\lambda \in \mln$, 
the holonomy map $\hol:P(S) \to R(S)$ 
is not injective on $\ov{\Q_\lam}$, 
although it is injective on $\Q_\lam$.  
\end{enumerate} 
\end{cor} 

Figure 2 is a computer generated figure by Komori, Sugawa, Wada, and 
Yamashita (cf. \cite{KoSu} and \cite{KoSuYaWa}), 
which is a part of a complex one-dimensional slice of $P(S)$ 
for a punctured-torus $S$. 
All the projective structures in this slice have the same 
underlying complex structure.  
The region in white is a slice of $Q(S)$: 
the inner disk is a Bers slice (a slice of the standard component $\Q_0$)  
and the outer part is a slice of an exotic component. 
This figure seems to illustrate the phenomenon 
stated in Theorem 1.3 and Corollary 1.4 (1).  
\begin{figure}[htbp]
  \begin{center}
    \includegraphics[keepaspectratio=true,height=44mm]{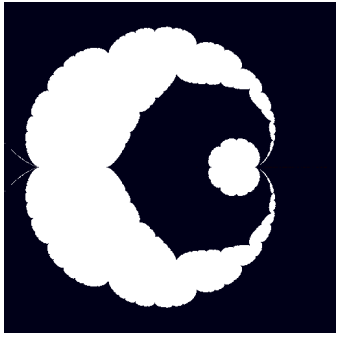}
  \end{center}
  \caption{A part of a slice of $Q(S)$ in $P(S)$
(white part) for a punctured-torus $S$.}
\end{figure}

We now consider associated results on the quasifuchsian space $\QF$ in $R(S)$. 
Let $\rho_n$ and $\rho_\fty$ denote holonomy representations 
of $Y_n$ and $Y_\fty$ respectively. 
Theorem 1.1 then implies that $\QF$ self-bumps at $\rho_\fty$ in $R(S)$; 
namely, for any neighborhood $U$ of $\rho_\fty$ in $R(S)$,  
there exits a neighborhood $V \subset U$ of $\rho_\fty$ such that 
$V \cap \QF$ is disconnected; see \cite{Mc}. 
As a direct consequence of Theorem 1.3 in the case of 
$i(\lam,\mu) \ne 0$, 
we can refine this statement as follows: 

\begin{thm}
For any neighborhood $U$ of $\rho_\fty$ in $R(S)$,  
there exits a neighborhood $V \subset U$ of $\rho_\fty$ 
such that for all large enough $n>0$, 
$\rho_n$ and $\rho_{-n}$ are contained in distinct components 
of $V \cap \QF$. 
\end{thm}

\begin{rem}
The same result is obtained independently by Bromberg and Holt 
(oral communication with Bromberg; see also Holt \cite{Ho}). 
\end{rem}

We mention here the relationship between the topics in this paper and 
the topology of deformation spaces of Kleinian groups.  
Let $\Gamma$ be a finitely generated Kleinian group with 
non-trivial space $AH(\Gamma)$ of conjugacy classes of 
discrete faithful representations $\rho:\Gamma \to \psl$.  
In general, the interior of $AH(\Gamma)$ 
consists of finitely many components. 
It was first shown by Anderson and Canary \cite{AnCa} 
that for some Kleinian group those components bumps. 
This result is generalized by Anderson, Canary and Metallurgy \cite{AnCaMc}. 
In our setting, 
$\Gamma$ is a quasifuchsian group isomorphic to $\pi_1(S)$, and 
the quasifuchsian space $\QF$ 
is the interior of the space $AH(\Gamma)$ 
which consists of exactly one connected component. 
Also in this case, 
the idea of Anderson and Canary can be applied to show that 
$\QF$ self-bumps by using projective structures; see Theorem 1.1 
due to McMullen. 
In fact, by lifting $\QF \subset R(S)$ to $Q(S) \subset P(S)$ 
via the holonomy map $\hol$, 
we can discuss the bumping 
and self-bumping of components of $Q(S)$.  
After McMullen, Bromberg and Holt \cite{BrHo} characterized 
the self-bumping of components of the interior 
of $AH(\Gamma)$ for general Kleinian groups $\Gamma$ 
without using projective structures. 
Our results, especially (1) and (2) in Corollary 1.4,  
can be viewed as the projective structure analogues of the works 
in \cite{AnCa}, \cite{AnCaMc} and \cite{BrHo}. 
We refer the reader to \cite{Ca} 
for further information on 
the bumping and self-bumping of deformation spaces of Kleinian groups.  

This paper is organized as follows: 
In section 2, we provide definitions and basic properties 
of the spaces and maps with which we are concerned. 
We devote section 3 to the proof of Theorem 1.3. 
The idea of the proof can be found at the beginning of this section. 
Corollaries of Theorem 1.3 are obtained in section 4. 

\section{Preliminaries}

We refer the reader to \cite{It1} for more information 
on some topics in this section. 
See also an exposition \cite{It2}. 

\subsection{Kleinian groups}

A {\it Kleinian group} $\Gamma$ is a discrete subgroup of $\psl$, 
which acts on the hyperbolic space ${\mathbb H}^3$ as isometries, and on 
the sphere at infinity $S^2_{\infty}=\wh{\C}$ 
as conformal automorphisms. 
The {\it region of discontinuity} 
$\Omega(\Gamma)$ is the largest open subset of 
$\wh{\C}$ on which $\Gamma$ acts properly discontinuously, 
and the {\it limit set} $\Lambda(\Gamma)$ of $\Gamma$ 
is its complement $\wh{\C}-\Omega(\Gamma)$. 
The quotient manifold 
$N_{\Gamma} =({\mathbb H}^3 \cup \Omega(\Gamma))/\Gamma$ is called 
the {\it Kleinian manifold} of $\Gamma$. 
A {\it quasifuchsian group} $\Gamma$ 
is a Kleinian group whose limit set $\Lambda(\Gamma)$ is a Jordan curve 
and which contains no element interchanging the two components 
of $\Omega(\Gamma)$. 
A {\it $b$-group} $\Gamma$ 
is a Kleinian group which has 
exactly one simply connected invariant component of $\Omega(\Gamma)$, 
which is denoted by $\Omega_0(\Gamma)$. 

Let $R(S)$ denote the space of all conjugacy classes $[\rho]$ of
representations of $\rho:\pi_1(S) \to \psl$ such that 
$\rho(\pi_1(S))$ is non-abelian. 
For simplicity, we denote $[\rho]$ by $\rho$ 
if there is no confusion. 
The space $R(S)$ is endowed with algebraic topology and  
is known to be a $(6g-6)$-dimensional complex manifold 
(see Theorem 4.21 in \cite{MaTa} for example). 
{\it Quasifuchsian space} $\QF=\QF(S)$ 
is the subset of $R(S)$ consisting of
faithful representations whose images are
quasifuchsian groups. 
Then it is known 
by Bers \cite{Be1}, Marden \cite{Ma} and Sullivan \cite{Su} that 
$\QF$ is connected, contractible and open in $R(S)$.  

\subsection{Space of projective structures}

A projective structure on $S$ is a $(G,X)$-structure 
where $X$ is the Riemann sphere $\widehat{\mathbb C}$ and 
$G=\psl$ is the group of projective automorphism of $\widehat{\mathbb C}$. 
Let $P(S)$ denote the space of marked projective structures on $S$, or 
the space of equivalence classes of pairs $(g, Y)$ 
of a projective surface $Y$ and an orientation 
preserving homeomorphism $g:S \to Y$.  
Two pairs $(g_1,Y_1)$ and $(g_2,Y_2)$ are said to be equivalent if 
there is an isomorphism $h:Y_1 \to Y_2$ of projective structures 
such that $h \circ g_1$ is isotopic to $g_2$. 
The equivalence class of $(g,Y)$ is simply denoted by $Y$.  

A projective structure $Y \in P(S)$ has an underlying 
conformal structure $\pi(Y) \in T(S)$, where 
$T(S)$ is the Teichm\"{u}ller space for $S$. 
The space $P(S)$ is equipped with a structure 
of a complex $(6g-6)$-dimensional 
holomorphic affine bundle over $T(S)$ 
with the projection $\pi:P(S) \to T(S)$. 

A projective structure $Y$ on $S$ 
can be lifted to that  
$\wt{Y}$ on $\wt{S}$, 
where $\wt{S} \to S$ is the universal cover on which 
$\pi_1(S)$ acts as a covering group.  
Since $\wt Y$ is simply connected, 
we obtain a developing map 
$f_Y:\wt Y \to \wh \C$ 
by continuing the charts analytically. 
In addition, the developing map induces a holonomy representation 
$\rho_Y: \pi_1(S) \cong \pi_1(Y) \to \psl$ 
which satisfies $f_Y \circ \gamma=\rho_Y(\gamma) \circ f_Y$ for 
every $\gamma \in \pi_1(S)$. 
We remark that the pair $(f_Y,\rho_Y)$ is determined 
uniquely up to the canonical action of $\psl$. 
We now define the {\it holonomy map} 
$$
\hol:P(S) \to R(S)
$$ 
by $Y \mapsto [\rho_Y]$. 
Then Hejhal \cite{He} showed that the map 
$\hol$ is a local homeomorphism and 
Earle \cite{Ea} and Hubbard \cite{Hu} independently
showed that the map is holomorphic:  

\begin{thm}[Hejhal, Earle and Hubbard]
The holonomy map $\hol:P(S) \to R(S)$ 
is a holomorphic local homeomorphism. 
\end{thm}

We denote by $Q(S)=\hol^{-1}(\QF)$ 
the set of projective structures with quasifuchsian holonomy. 
An element of $Q(S)$ is said to be {\it standard} if 
its developing map is injective; otherwise it is {\it exotic}. 
Let $\Q_0 \subset Q(S)$ denote 
the set of all standard projective structures. 
Then the map $\hol|_{\Q_0}:\Q_0 \to \QF$ is a biholomorphism, 
which takes the Bers' embedded image of the Theichm\"{u}ller space $T(S)$ 
in a fiber to a Bers slice (cf. \cite{Be2}). 

\subsection{Integral points of measured laminations} 

Let $\mathcal S=\mathcal S (S)$ 
denote the set of homotopy classes of non-trivial 
simple closed curves on $S$.  
By abuse of the notation, we also denote a representative 
of $c \in {\mathcal S}$ by $c$. 
Let $\mln=\mln(S)$ denote the set of 
integral points of measured laminations on $S$, 
or the set of formal summation $\sum_{i=1}^l k_ic_i$ of mutually 
distinct, disjoint elements $c_i \in \mathcal S$  
with positive integer $k_i$ weights. 
We regard $\mathcal S \subset \mln$. 
The ``zero-lamination'' $0$ is also contained in $\mln$. 
A {\it realization} $\widehat{\lam}$ of 
$\lam = \sum_{i=1}^l k_ic_i \in \mln$ 
is a disjoint union of simple closed curves 
which realize each weighted simple closed curve $k_ic_i$ 
by $k_i$-parallel disjoint simple closed curves homotopic to $c_i$. 

For $c,\,d \in {\mathcal S}$, the geometric intersection number 
$i(c,d)$ is the minimum number of points in which 
the representations of $c$ and $d$ must intersect. 
Note that $i(c,c)=0$ for any $c \in {\mathcal S}$. 
We naturally extend the definition of the geometric intersection number 
for elements of $\mln$ by linearity. 
If $i(\lambda, \mu) =0$ for $\lam,\,\mu \in \mln$, 
we can define $m \lambda + n \mu \in \mln$ 
for $m,n \in \N$. 

\subsection{Grafting}

Let $Y \in P(S)$ and $\wtil{Y} \to Y$ the universal cover. 
A simple closed curve $c \in {\mathcal S}$ 
is called {\it admissible} on $Y$ 
if a connected component $\wt{c} \subset \wt Y$ 
of the preimage of $c \subset Y$ 
is mapped injectively by the developing map $f_Y:\wt Y \to \wh\C$ 
onto its image $f_Y(\wt c)$ 
and if the holonomy image $\rho_Y(c)$ fixing $f_Y(\wt c)$ 
is loxodromic. 
We say that $\lam \in \mln$ is admissible on $Y$ if 
every component of the support of $\lam$ is admissible.  

Suppose that $c \in \cal S$ is admissible on $Y \in P(S)$. 
Let 
$$
A_c=(\wh \C -f_Y(\wt c))/\langle \rho_Y(c) \rangle 
$$
be the quotient annulus with its induced projective structure.  
Then the {\it grafting} $\Gr_c(Y)$ is the projective surface 
obtained by cutting $Y$ along $c$ and inserting 
the annulus $A_c$ without twisting. 
Similarly, we can define the grafting 
$\Gr_\lam(Y)$ for admissible $\lam \in \mln$ by linearity.   
The basic fact is that the grafting operation does not change the holonomy 
representation; i.e. $\hol(\Gr_\lam(Y))=\hol(Y)$ holds for every 
admissible $\lam \in \mln$.  
Since the map $\hol|_{\Q_0}:\Q_0 \to \QF$ is a biholomorphism 
and the grafting map 
$$
\Gr_\lam:\Q_0 \to P(S)
$$
for $\lam \in \mln$ satisfies $hol \circ \Gr_\lam=hol$,  
the map $\Gr_\lam$ turns out to be a biholomorphism 
from $\Q_0$ onto its image $\Q_\lam=\Gr_\lam(\Q_0)$. 
From Goldman's grafting theorem \cite{Go} below, 
we obtain the decomposition of $Q(S)$ into its 
connected components; 
$$
Q(S)=\bigsqcup_{\lam \in \mln}
\cal Q_\lam. 
$$   

\begin{thm}[Goldman]
Suppose that $Y \in \Q_0$. 
Then every projective structure with holonomy $\hol(Y)$ 
is obtained by grafting of $Y$ along some $\lam \in \mln$. 
\end{thm}

\subsection{Operations on $\mln$}

We now introduce two operations 
$$
(\cdot,\cdot)_{\sharp},\, (\cdot,\cdot)_{\flat}: \mln \times \mln \to \mln,
$$ 
which is closely observed by Luo in \cite{Lu}. 
For $\lam,\,\mu \in \mln$, 
the elements $(\lam,\mu)_{\sharp}$ and 
$(\lam,\mu)_{\flat}$ in $\mln$ are obtained as follows: 
Let $\what{\lam}, \,\what{\mu}$ be realizations of $\lam, \, \mu$ whose 
geometric intersection number is minimal. 
Next resolve all intersecting points in  $\what{\lam} \cup\what{\mu}$ 
as in Figure 3.  
Then $(\lam,\mu)_{\sharp}$ and 
$(\lam,\mu)_{\flat}$ are elements in $\mln$ whose realizations 
are the resulting curve systems, respectively; see also Figure 4.
We collect in the next lemma some basic properties of these operations, 
whose proof is left for the reader. 

\begin{lem}
For $\lam, \mu \in \mln$, we have:  
\begin{enumerate}
\item 
$(\lam,\mu)_{\sharp}=(\mu,\lam)_{\flat}$. 
\item 
$(\lam,\mu)_{\sharp} = (\lam,\mu)_{\flat}$ if and only if $i(\lam,\mu) = 0$.  
If $i(\lam,\mu)=0$ then  $(\lam,\mu)_{\sharp}=(\lam,\mu)_{\flat}=\lam+\mu$. 
\item 
Suppose that every component of $\mu$ intersects $\lam$ essentially. 
Then 
$((\lambda,\mu)_{\sharp}, \mu)_{\flat}=((\lambda,\mu)_{\flat}, \mu)_{\sharp}=\lambda$. 
\end{enumerate}
\end{lem}
\begin{figure}[htbp]
\begin{center}
\includegraphics{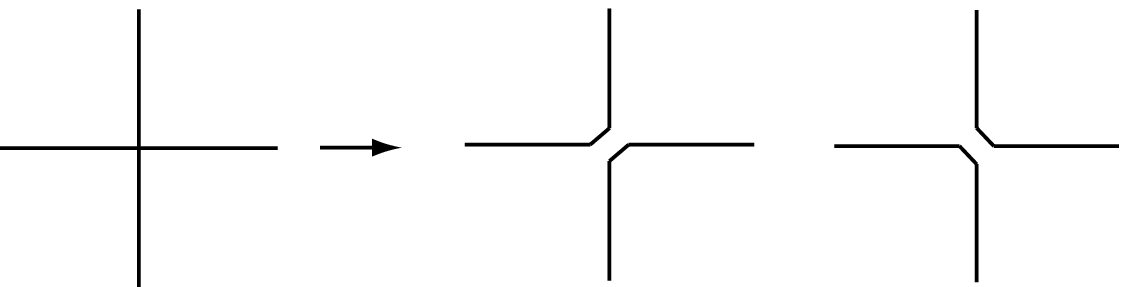}
\end{center}
\vspace{-4.5cm}
\begin{center}
\unitlength 0.1in
\begin{picture}( 57.5000, 16.5000)(  6.0000,-21.5000)
%
\special{pn 8}%
\special{pa 600 500}%
\special{pa 6350 500}%
\special{pa 6350 2150}%
\special{pa 600 2150}%
\special{pa 600 500}%
\special{ip}%
\put(10.8000,-10.6000){\makebox(0,0)[lb]{$\wh \lam$}}%
\put(18.1000,-18.6000){\makebox(0,0)[lb]{$\wh \mu$}}%
\put(33.5000,-20.8000){\makebox(0,0)[lb]{$\wh{(\lam,\mu)_\sharp}$}}%
\put(48.3000,-20.7000){\makebox(0,0)[lb]{$\wh{(\lam,\mu)_\flat}$}}%
\end{picture}%
\end{center}
\caption{Rules to obtain $(\lam,\mu)_{\sharp}$ 
and $(\lam,\mu)_{\flat}$.}
\end{figure}
\begin{figure}[htbp]
\begin{center}
\includegraphics{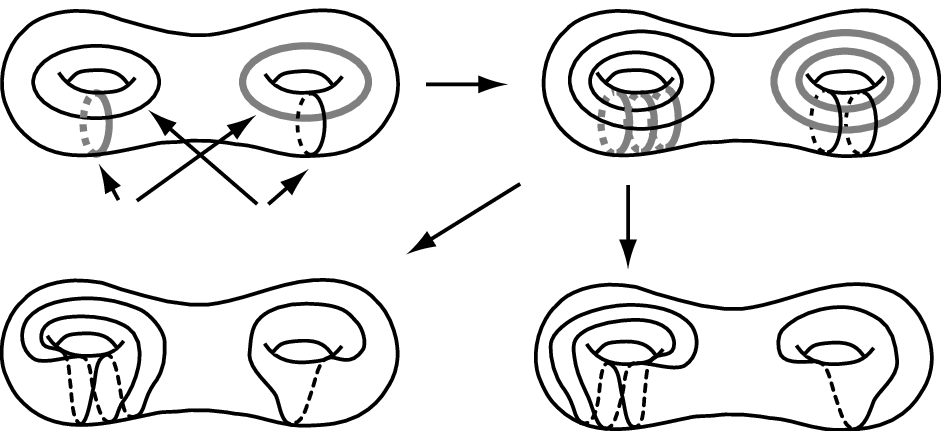}
\end{center}
\vspace{-6cm}
\begin{center}
\unitlength 0.1in
\begin{picture}( 49.6000, 23.3000)(  6.3000,-28.2000)
%
\special{pn 8}%
\special{pa 630 490}%
\special{pa 5590 490}%
\special{pa 5590 2820}%
\special{pa 630 2820}%
\special{pa 630 490}%
\special{ip}%
\put(18.7000,-10.4000){\makebox(0,0)[lb]{$2$}}%
\put(21.6000,-10.4000){\makebox(0,0)[lb]{$2$}}%
\put(17.4000,-13.6000){\makebox(0,0)[lb]{$3$}}%
\put(25.6000,-13.8000){\makebox(0,0)[lb]{$2$}}%
\put(16.7000,-17.5000){\makebox(0,0)[lb]{$\lam$}}%
\put(22.1000,-17.5000){\makebox(0,0)[lb]{$\mu$}}%
\put(39.7000,-16.2000){\makebox(0,0)[lb]{$\wh \lam$ and $\wh \mu$}}%
\put(17.7000,-27.4000){\makebox(0,0)[lb]{$(\lam,\mu)_\sharp$}}%
\put(39.9000,-27.4000){\makebox(0,0)[lb]{$(\lam,\mu)_\flat$}}%
\put(23.5000,-23.7000){\makebox(0,0)[lb]{$2$}}%
\put(19.1000,-23.5000){\makebox(0,0)[lb]{$1$}}%
\put(40.2000,-23.9000){\makebox(0,0)[lb]{$1$}}%
\put(46.1000,-23.8000){\makebox(0,0)[lb]{$2$}}%
\end{picture}%
\end{center}
\caption{Examples of $(\lam,\mu)_{\sharp}$ and $(\lam,\mu)_{\flat}$.}
\end{figure}

\subsection{Geometric limits of Kleinian groups}

We begin with the definition of Hausdorff convergence. 

\begin{defn}[Hausdorff convergence, geometric convergence] 
Let $X$ be a locally compact Hausdorff space. 
A sequence of closed subsets $A_n \subset X$ is said to converge 
in $X$ to 
a closed subset $A \subset X$ 
{\it in the sense of Hausdorff} 
if every element $x \in A$ is the limit of a sequence 
$\{x_{n} \in A_n \}$ and if every accumulation point of every 
sequence $\{x_{n} \in A_n \}$ lies in $A$. 
This is also denoted by
$A_n \hto A$ in $X$. 
A sequence of Kleinian groups $\Gamma_n$ is said to converge 
{\it geometrically} to a group $\wh{\Gamma}$ if 
$\Gamma_n$ converges in $\psl$ to $\wh \Gamma$ in the sense of Hausdorff. 
\end{defn} 

It is a result of J{\o}rgensen and Marden \cite{JoMa} that 
if $\rho_n \to \rho_\fty$ in $AH(S)$ then 
the sequence $\Gamma_n=\rho_n(\pi_1(S))$ 
converges geometrically to a Kleinian group $\wh{\Gamma}$ 
up to taking a subsequence. Furthermore, 
the geometric limit $\wh{\Gamma}$ always 
contains the algebraic limit $\Gamma_{\infty}=\rho_\infty(\pi_1(S))$. 
The following theorem is due to Kerckhoff and Thurston \cite{KeTh}. 

\begin{thm}[Kerckhoff-Thurston] 
Suppose that a sequence $\rho_n \in \QF$ 
converges to some $\rho_\infty \in AH(S)$ 
and that the sequence $\Gamma_n=\rho_n(\pi_1(S))$ 
converges geometrically to $\wh{\Gamma}$. 
Then the sequence $\Lam(\Gamma_n)$ converges in $\wh{\C}$ 
to $\Lam(\wh{\Gamma})$ in the sense of Hausdorff. 
\end{thm} 

\subsection{Pullbacks of limit sets of Kleinian groups} 

Let $Y \in P(S)$ be a projective structure 
with discrete faithful holonomy 
$\rho_Y:\pi_1(S) \to \Gamma$.  
Let $\pi_Y:\wt{Y} \to Y$ be the universal cover and let 
$f_Y:\wt{Y} \to \wh{\C}$ be the developing map. 
Then the preimage $f_Y^{-1}(\Lam(\Gamma))$  of the limit set 
$\Lam(\Gamma)$ in $\wt{Y}$ is invariant under the action of the 
covering transformation group $\pi_1(Y)$.  
Thus the subset $f_Y^{-1}(\Lam(\Gamma))$ in $\wt{Y}$ descends to 
the subset
$$
\Lam_Y:=\pi_Y \circ f_Y^{-1}(\Lam(\Gamma))  
$$
in $Y$, which is called the {\it pullback} of the limit set 
$\Lam(\Gamma)$ in $Y$. 
Similarly, we also   
obtain the pullback $\pi_Y \circ f_Y^{-1}(\Lam(\wh{\Gamma}))$ 
of the limit set $\Lam(\wh{\Gamma})$ 
for any Kleinian group $\wh{\Gamma}$ containing $\Gamma$. 
By definition of grafting and Goldman's grafting theorem, 
an element $Y \in \Q_\lam$ is characterized as follows: 

\begin{lem}
Let $Y$ be an element of $Q(S)$. 
Then $Y \in \Q_{\lambda}$ if and only if  $\Lambda_Y \subset Y$ is a 
realization of $2\lambda$.  
\end{lem}

Suppose that $Y_n \to Y_\fty$ in $P(S)$ as $n \to \infty$. 
A projective structure on $S$ induces a complex structure on $S$,  
and hence a hyperbolic structure. 
Therefore, with these canonical hyperbolic structures on $Y_n$ and $Y_{\infty}$, 
there exist $K_n$-quasi-isometric maps 
$\omega_n:Y_{\infty} \to Y_n$ with $K_n \to 1$ as $n \to \infty$. 

\begin{lem}[Lemma 3.3 in {\cite{It1}}]
Suppose that a sequence $Y_n \in Q(S)$ 
converges to some $Y_\fty \in \ov{Q(S)}$ as $n \to \infty$. 
Then $\rho_{Y_n} \to \rho_{Y_\fty}$ in $R(S)$. 
We further assume that the sequence 
$\Gamma_n=\rho_{Y_n}(\pi_1(S))$ converges geometrically 
to a Kleinian group $\wh{\Gamma}$. 
Now let $\omega_n:Y_\fty \to Y_n$ be a $K_n$-quasi-isometric map with  
$K_n \to 1$ as $n \to \infty$. 
Then the sequence 
$\omega_{n}^{-1}(\Lam_{Y_n})$ converges in $Y_\fty$ 
to $\wh{\Lam}_{Y_\fty}$ 
in the sense of Hausdorff, where 
$\wh{\Lam}_{Y_\fty}=\pi_{Y_\fty} \circ f_{Y_\fty}^{-1}(\Lam(\wh{\Gamma}))$ 
is the pullback of the limit set $\Lam(\wh \Gamma)$ in $Y_\fty$. 
If there is no confusion, we simply say that 
$\Lambda_{Y_n}$ converge in $Y_\fty$ to $\wh{\Lam}_{Y_\fty}$. 
\end{lem} 

\section{Proof of the main theorem}

We devote this section to the proof of Theorem 1.3. 
Throughout this proof, Figure 5 should be helpful for the reader to 
understand the arguments.  
Here we outline the proof.  
Let $\rho_n,\,Y_n$ and $Z_n$ be sequences as in introduction.  
We will recall in \S 3.1--3.2 constructions and basic facts 
of sequences $\rho_n$ and $Y_n$; 
see \cite{Mc} and \cite{It1} for more details. 
Suppose that the sequence $\Gamma_n=\rho_n(\pi_1(S))$ 
converges geometrically to a Kleinian group $\wh \Gamma$. 
To show that $\{Z_n\}_{n \gg 0} \subset \Q_{(\lam,\mu)_\sharp}$ 
and $\{Z_n\}_{n \ll 0} \subset \Q_{(\lam,\mu)_\flat}$, by Lemma 2.6, 
it suffices to show that 
the pullback $\Lam_{Z_n}$ of 
$\Lam(\Gamma_n)$ in $Z_n$ is a realization of 
$2(\lam,\mu)_\sharp$ for all $n \gg 0$ and 
of $2(\lam,\mu)_\flat$ for all $n \ll 0$.   
Note that $\Lam_{Z_n} \subset Z_n$ converge 
to the pullback 
$\LZ \subset Z_\fty$ of $\Lam(\wh \Gamma)$ in the sense of Lemma 2.7. 
Thus to understand the shape of $\Lam_{Z_n}$ 
we first study in \S 3.4 the shape of $\LZ \subset Z_\fty$ 
by using the result for the shape of 
the pullback $\LY \subset Y_\fty$ of $\Lam(\wh \Gamma)$ 
obtained in \cite{It1} (see \S 3.3).  
We will see in \S 3.5 that 
$\Lam_{Z_n}$ are obtained by modifying $\LZ$ at some points 
which are pullbacks of rank-two parabolic fixed points in $\Lam(\wh \Gamma)$. 
The reason for the difference of the resulting curve systems 
$\Lam_{Z_n}$ for $n \gg 0$ 
and for $n \ll 0$ is that the limit sets 
$\Lam(\Gamma_n)$ for $n \gg 0$ and for $n \ll 0$ are spiraling in opposite directions 
at each rank-two parabolic fixed point in $\Lam(\wh \Gamma)$, 
which will be studied closely in \S 3.5. 
\begin{figure}[htbp]
\begin{center}
\includegraphics{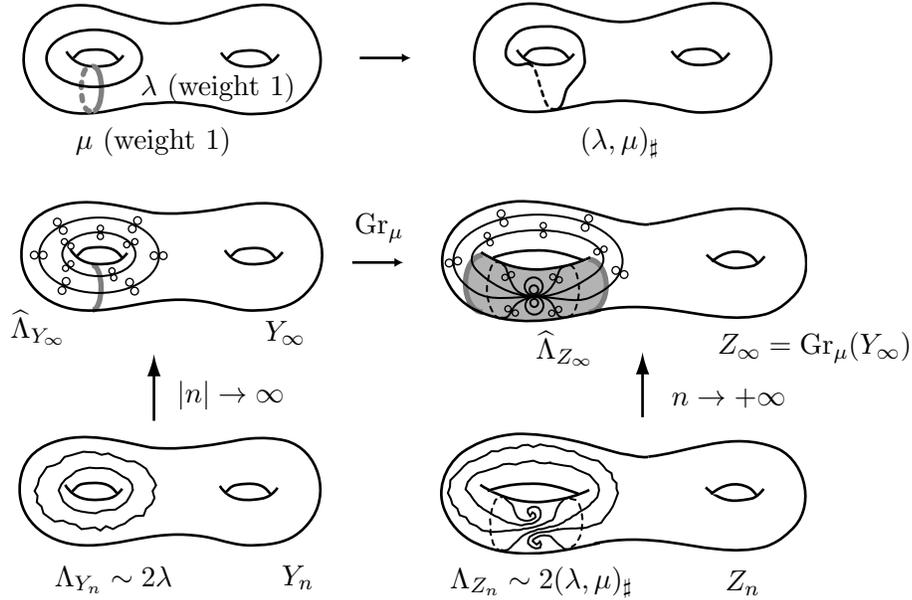}
\end{center}
\vspace{-9cm}
\begin{center}
\unitlength 0.1in
\begin{picture}( 50.5000, 35.7000)(  7.6000,-40.7000)
%
\special{pn 8}%
\special{pa 760 500}%
\special{pa 5810 500}%
\special{pa 5810 4070}%
\special{pa 760 4070}%
\special{pa 760 500}%
\special{ip}%
\put(18.6000,-13.3000){\makebox(0,0)[lb]{$\lambda$ (weight $1$)}}%
\put(15.2000,-16.1000){\makebox(0,0)[lb]{$\mu$ (weight $1$)}}%
\put(41.6000,-16.2000){\makebox(0,0)[lb]{$(\lam,\mu)_{\sharp}$}}%
\put(29.8000,-20.6000){\makebox(0,0)[lb]{$\Gr_{\mu}$}}%
\put(25.1000,-26.0000){\makebox(0,0)[lb]{$Y_\fty$}}%
\put(48.8000,-26.9000){\makebox(0,0)[lb]{$Z_\fty=\Gr_{\mu}(Y_\fty)$}}%
\put(26.0000,-38.8000){\makebox(0,0)[lb]{$Y_n$}}%
\put(49.2000,-39.1000){\makebox(0,0)[lb]{$Z_n$}}%
\put(11.8000,-26.0000){\makebox(0,0)[lb]{$\LY$}}%
\put(39.3000,-27.1000){\makebox(0,0)[lb]{$\LZ$}}%
\put(14.1000,-39.0000){\makebox(0,0)[lb]{$\Lambda_{Y_n} \sim {2 \lambda}$}}%
\put(34.8000,-39.2000){\makebox(0,0)[lb]{$\Lambda_{Z_n} \sim {2 (\lambda,\mu)_{\sharp}}$}}%
\put(20.5000,-29.4000){\makebox(0,0)[lb]{$|n| \to \infty$}}%
\put(46.4000,-29.3000){\makebox(0,0)[lb]{$n \to +\fty$}}%
\end{picture}%
\end{center}
\caption{Schematic figure explaining the proof of Theorem 1.3.}
\end{figure}

\subsection{A Kleinian group with rank-two parabolic subgroups}

Throughout this section, we fix a non-zero element 
$\lam=\sum_{i=1}^l k_i c_i \in \mln$ with its support 
$\ul{\lam}=\sqcup_{i=1}^l c_i$. 
Let 
$$
M_\lam=S \times [-1,1] - \ul{\lam} 
\times \left\{0\right\}
$$ 
be a 3-manifold $S \times [-1,1]$ with simple closed curves 
$c_i \times \{0\} \, (1 \le i \le l)$ removed. 
Let $\wh{\Gamma}$ be a Kleinian group 
whose Kleinian manifold 
$N_{\wh{\Gamma}}=({\mathbb H}^3 \cup \Omega(\wh{\Gamma}))/\wh{\Gamma}$ 
is homeomorphic to $M_\lam$. 
In what follows,  
we identify $N_{\wh{\Gamma}}$ with $M_\lam$ via this homeomorphism. 
Each tubular neighborhood of $c_i \times \{ 0 \}$ in $M_\lam$ 
corresponds to a rank-two cusp end in $N_{\wh \Gamma}$,  
and hence to a conjugacy class of  maximal 
rank-two parabolic subgroup 
$\langle \gamma_i,\delta_i \rangle$ in 
$\wh{\Gamma} \cong \pi_1(N_{\wh \Gamma})$. 
We fix the generators of the group   
$\langle \gamma_i,\delta_i \rangle$ 
so that $\gamma_i \in \wh \Gamma$ is freely homotopic to 
$c_i \times \{-1\}$ in $N_{\wh \Gamma}$ 
and that $\delta_i \in \wh \Gamma$ is trivial in $S \times [-1,1]$. 
Moreover, we orient $\gamma_i$ and $\delta_i$ 
so that for some $\phi_i \in \psl$, 
$\phi_i \circ \gamma_i \circ \phi_i^{-1}(z)=z+1$ 
and 
$\phi_i \circ \delta_i \circ \phi_i^{-1}(z)=z+\tau_i$ 
with $\Im \tau_i>0$. 

Note that each connected component $\omega$ of $\Omega(\wh \Gamma)$ 
covers a connected component of the conformal boundary 
$\bd N_{\wh \Gamma}=S \times \{\pm 1\}$ via the 
quotient map 
${\mathbb H}^3 \cup \Omega(\wh{\Gamma}) \to N_{\wh{\Gamma}}$, and that 
the subgroup $\Gamma$ of $\wh \Gamma$ 
stabilizing $\omega$ is a $b$-group with $\omega=\Omega_0(\Gamma)$. 

\subsection{Wrapping maps and associated representations}

We introduce here the wrapping map $w_\lam:S \to N_{\wh{\Gamma}}$ 
associated to $\lam \in \mln$, which is 
an immersion determined up to homotopy. 
For $0 \in \mln$, we let 
$w_0:S \to S \times\{-1/2\} \subset N_{\wh \Gamma}$ be the 
canonical inclusion. 
The wrapping map $w_\lam:S \to N_{\wh \Gamma}$ associated to 
$\lam \in \sum_{i=1}^l k_i c_i \in \mln$ 
is an immersion such that 
the image $w_\lam(S)$ in $N_{\wh \Gamma}$ is obtained by cutting 
$S \times \{-1/2\}$ along $c_i \times \{-1/2\}$ and inserting an annulus 
which wraps $k_i$-times around $c_i \times \{0\}$ in $N_{\wh \Gamma}$ 
at the cut locus for every $1 \le i \le l$; see Figure 6. 
It is also required that $w_\lam$ is homotopic to $w_0$ in $S \times [-1,1]$.  
\begin{figure}[htbp]
\begin{center}
\includegraphics{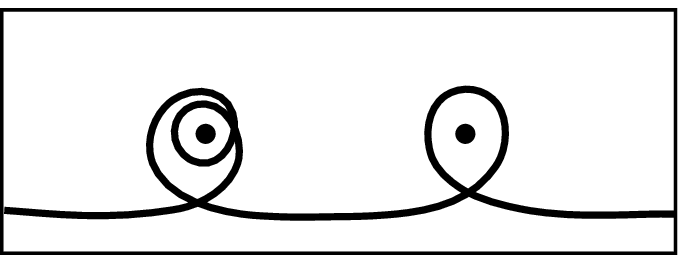}
\end{center}
\vspace{-3.5cm}
\begin{center}
\unitlength 0.1in
\begin{picture}( 45.9000, 11.5000)( 17.1000,-25.2000)
%
\special{pn 8}%
\special{pa 1710 1370}%
\special{pa 6300 1370}%
\special{pa 6300 2520}%
\special{pa 1710 2520}%
\special{pa 1710 1370}%
\special{ip}%
\put(54.7000,-15.5000){\makebox(0,0)[lb]{$S \times \{1\}$}}%
\put(54.7000,-24.6000){\makebox(0,0)[lb]{$S \times \{-1\}$}}%
\put(54.6000,-20.0000){\makebox(0,0)[lb]{$N_{\what{\Gamma}} \subset S \times [-1,1]$}}%
\put(31.0000,-16.5000){\makebox(0,0)[lb]{$c_i \times \{0\}$}}%
\put(43.5000,-16.5000){\makebox(0,0)[lb]{$c_j \times \{0\}$}}%
\put(37.5000,-21.8000){\makebox(0,0)[lb]{$w_\lam(S)$}}%
\end{picture}%
\end{center}
\caption{Schematic figure of a wrapping map $w_\lam:S \to N_{\what{\Gamma}}$. }
\end{figure}
Here we make use of Thurston's Dehn filling theorem 
(\cite{Co}, see also \cite{BoOt}).  
By performing simultaneous $(-1,n)$ Dehn filling ($n \in \Z$) 
on every cusp end of $N_{\wh{\Gamma}}$, we obtain a 
sequence of representations $\{\chi_n:\wh{\Gamma} \to \psl\}$  
which satisfies the following: 
\begin{itemize} 
\item $\Gamma_n=\chi_n(\wh{\Gamma})$ is a quasifuchsian group, 
\item The kernel of $\chi_n$ is normally generated by 
$\gamma_1^{n} \delta_1^{-1}, \ldots ,\gamma_l^{n} \delta_l^{-1}$, 
\item $\chi_n$ converge algebraically to 
the identity map of $\wh \Gamma$ as $|n| \to \infty$, and 
\item $\Gamma_n$ converge geometrically to $\wh{\Gamma}$ as $|n| \to \infty$. 
\end{itemize} 
Then $\chi_n(\gamma_i) \to \gamma_i$ 
and ${\chi_n(\gamma_i)}^{n}=\chi_n({\gamma_i}^n)
=\chi_n(\delta_i) \to \delta_i$ 
as $|n| \to \infty$ for each $1 \le i \le l$. 
Now set  
\begin{eqnarray*}
\rho_n=\chi_n \circ (w_{\lambda})_*&:\pi_1(S) \to \psl, \\
\rho_{\infty}=(w_{\lambda})_*&:\pi_1(S) \to \psl, 
\end{eqnarray*}
where $(w_{\lambda})_*:\pi_1(S) \to \pi_1(N_{\wh{\Gamma}})=\wh{\Gamma}$ 
is the group isomorphism induced by $w_{\lambda}$. 
Then $\rho_n$ are faithful representations onto the quasifuchsian groups  
$\Gamma_n$ and 
the sequence $\rho_n$ converges algebraically to $\rho_{\infty}$. 
The algebraic limit $\Gamma_\fty=\rho_{\infty}(\pi_1(S))$ 
is a proper subgroup of the geometric limit $\wh{\Gamma}$, 
which is a geometrically finite $b$-group 
whose Kleinian manifold $N_{\Gamma_\fty}$ is homeomorphic to 
$S \times [-1,1]-\ul{\lam} \times \{1\}$. 
Let 
$$
Y_\fty=\Omega_0(\Gamma_\fty)/\Gamma_\fty
$$ 
be the projective structure 
on the conformal end of $N_{\Gamma_\fty}$ corresponding to 
$S \times \{-1\}$. 
Then $Y_\fty  \in \bd \Q_0$ and $\hol(Y_\fty)=\rho_\fty$. 
Since the map $\hol$ is a local homeomorphism, 
there is a sequence 
$\{Y_n\}_{|n| \gg 0}$ in $Q(S)$ which satisfies 
$Y_n \to Y_\fty$ as $|n| \to \fty$ 
and $\hol(Y_n)=\rho_n$ for all $|n| \gg 0$. 
Then it is known by McMullen \cite{Mc} 
that $Y_n$ are exotic for all $|n| \gg 0$; see Theorem 1.1. 

\subsection{The pullback of the limit set in $Y_\fty$} 

Let $\pi_{Y_\fty}:\wt Y_\fty=\Omega_0(\Gamma_\fty) \to Y_\fty$ be the 
universal cover. 
We recall from \cite{It1} the shape of the pullback 
$$
\LY=\pi_{Y_\fty} \circ f_{Y_\fty}^{-1}(\Lam(\wh \Gamma))=
(\Lam(\wh \Gamma) \cap \Omega_0(\Gamma_\fty))/\Gamma_\fty
$$
of the limit set $\Lam(\wh \Gamma)$ in $Y_\fty$. 
We first fix our terminology: 
A {\it regular neighborhood} $\rn (\xi)$ of 
a finite union of (not necessarily disjoint) 
simple closed curves  $\xi$ on $S$ 
is a 2-dimensional submanifold of $S$ 
with a deformation retraction 
$r:\rn (\xi) \to \xi$. 

\begin{prop}[Lemma 4.1 in \cite{It1}]
There is a regular neighborhood 
$\cal N(\what{\lam})$ of a realization 
$\wh \lam$ of $\lam$ in $Y_\infty$
which satisfy the following: 
\begin{enumerate}
\item $\LY$ is contained in the interior of $\rn (\wh{\lam})$, and 
\item for each connected component $\rn$ of $\cal N(\what{\lam})$, 
$\rn \cap \LY$ consists of two essential 
connected component. 
\end{enumerate}
\end{prop}

Since $\Lam_{Y_n}$ converge  in $Y_\fty$ to $\LY$ in the sense of Lemma 2.7, 
we can deduce from Proposition 3.1 
that $\Lam_{Y_n}$ are realization of $2\lam$, 
and hence from Lemma 2.6 that 
$Y_n \in \Q_\lam$ for all large $|n|$, 
which implies Theorem 1.2; see \cite{It1} for more details. 

In what follows, 
we will also make use of a regular neighborhood 
$\cal N(\underline{\lam})$ of the support $\ul{\lam}$ of $\lam$ 
which contains $\cal N(\what{\lam})$, 
where a connected component $\rn (c_i)$ of $\rn (\ul{\lam})$ 
contains $k_i$-parallel components of $\rn(\wh{\lam})$.  
Throughout this proof, 
we fix these regular neighborhoods: 
$$
\LY \subset \rn (\wh{\lam}) \subset \rn (\ul{\lam}) \subset Y_\fty. 
$$
Note that Proposition 3.1 (1) implies that 
$\Lambda(\wh{\Gamma}) \cap \Omega_0(\Gamma_\fty)$ 
lies in $\pi_{Y_\fty}^{-1}(\rn (\wh{\lambda}))$, 
and hence in $\pi_{Y_\fty}^{-1}(\rn (\ul{\lambda}))$. 

\subsection{The pullback of the limit set in $Z_\fty$}

Let $\mu$ be a non-zero element of $\mln$ 
which has no parallel component in common with $\lam$. 
Then $\mu$ is admissible on $Y_\fty$, and thus 
the grafting  $Z_\fty=\Gr_\mu(Y_\fty)$ of $Y_\fty$ along $\mu$ 
is obtained. 
We remark that $Z_\fty$ lies in $\bd \Q_{\mu}$ 
because the map 
$\Gr_\mu:\Q_0 \to \Q_\mu$ extends continuously to some 
neighborhood of $Y_\fty \in \bd \Q_0$; see \cite{Ita}. 
Let $\pi_{Z_\fty}:\wt Z_\fty \to Z_\fty$ be the 
universal cover. 
We now study the shape of the pullback 
$$
\LZ=\pi_{Z_\fty} \circ f_{Z_\fty}^{-1}(\Lam(\wh \Gamma))
$$
of $\Lam(\wh \Gamma)$ in $Z_\fty$ 
by using the observation of the shape of $\LY \subset Y_\fty$ in \S 3.3. 

\begin{prop}
There is a regular neighborhood ${\cal N}(\wh{\lam} \cup \wh{\mu})$ 
of $\wh \lam \cup \wh \mu$ in $Z_\fty$ 
which contains $\LZ$, 
where $\wh \lam,\, \wh \mu$ are realizations of $\lam,\,\mu$ 
whose geometric intersection number is minimal. 
\end{prop}

We devote this subsection to the proof of this proposition. 
For simplicity, we will show in the case where $\mu$ 
is a simple closed curve $d \in \cal S$ of weight one. 
Recall that $\pi_{Y_\fty}:\wt Y_\fty=\Omega_0(\Gamma_\fty) \to Y_\fty$ 
is the universal cover. 
Since $d \subset Y_\fty$ is admissible, 
a connected component $\wt{d}$ of  
$\pi_{Y_\fty}^{-1}(d)$ in $\Omega_0(\Gamma_\fty)$ is 
$\langle \eta \rangle$-invariant for some  
loxodromic element $\eta \in \Gamma_\fty$. 
Let 
$$
T=(\wh \C -\fix(\eta))/\langle \eta \rangle
$$ 
be the quotient torus and let 
$$
\pi_T:\wh \C -\fix(\eta) \to T
$$ 
be the covering map. 
The simple closed curve 
$\pi_T(\wt d)$ in $T$ is also denoted by $d$. 
We set the longitude $l$ of $T$ 
as a unique simple closed curve on $T$ 
(up to homotopy) which is contractible in 
${\mathbb H}^3/ \langle \eta \rangle$, 
while $d$ is the meridian of $T$. 

Recall that the grafting $Z_\fty=\Gr_d(Y_\fty)$ is obtained by 
cutting $Y_\fty$ along $d$ and inserting a cylinder $T-d$:  
\begin{eqnarray*}
Z_\fty=(Y_\fty - d) \sqcup (T - d). 
\end{eqnarray*}
Then the developing map 
$f_{Z_\fty}:\wt{Z}_\fty \to \hatC$ is obtained by 
gluing piece by piece the developing maps 
\begin{eqnarray*}
f_{Y_\fty}:\wt{Y}_\fty=\Omega_0(\Gamma_\fty) \to \hatC,  \quad 
f_{T}: \wt T=\hatC-\fix (\eta) \to \hatC 
\end{eqnarray*}
and  their conjugations by elements of $\Gamma_\fty$. 
(Note that both $f_{Y_\fty}$ and $f_{T}$ are equal to the identity map. )
Therefore the set 
$\LZ \cap (Y_\fty - d)$  
is equal to $\LY \setminus d$ and is contained in 
$\rn (\wh \lam) \setminus d$ as a subset of $Y_\fty-d \subset Z_\fty$, 
and the set 
$\LZ \cap (T - d)$ is equal to 
$\wh \Lam_{T} \setminus d$ 
as a subset of $T-d \subset Z_\fty$, where 
$$
\wh \Lam_{T}=\pi_T \circ 
f_{T}^{-1}(\Lam(\wh \Gamma))=
\pi_T(\Lam(\wh \Gamma)-\fix (\eta)) 
$$
is the pullback of $\Lam(\wh \Gamma)$ in $T$. 
Thus we have 
\begin{eqnarray*}
\LZ=(\LY \setminus d) \sqcup (\wh \Lam_{T} \setminus d). 
\end{eqnarray*}

We first consider the case of $i(\lam,d)=0$. 
We may assume that $\rn (\wh{\lam}) \cap d=\emptyset$.   
Then we have $\LY \cap d=\emptyset$ in $Y_\fty$ 
from $\LY \subset \rn (\wh{\lam})$. 
Moreover, this implies that 
$\Lam(\wh{\Gamma}) \cap \wt{d}=\emptyset$ in $\hatC$, and hence that 
$\wh \Lam_{T} \cap d=\emptyset$ in $T$. 
Thus we obtain  
$\LZ \subset \rn (\wh{\lam} \cup d)$ by setting 
a regular neighborhood $\rn (\wh{\lam} \cup d)$ of 
$\wh{\lam} \cup d$ as 
$$
\rn (\wh{\lam} \cup d)= \rn (\wh{\lam}) \sqcup \rn_{T}, 
$$
where $\rn_{T}$ is an annulus in $T-d \subset Z_\fty$ 
which contains $\LT$. 

From now on, we assume that $i(\lam, \mu)=i(\lam,d) \ne 0$ 
and set $s=i(\ul{\lam}, d)$. 
Here we fix our terminology: 
\begin{defn}[crescent]
A closed set $A \subset \wh \C$ is called a {\it crescent}  
if $A$ is the closure of $B_2- B_1$,  
where  $B_1, \,B_2$ are 
topological closed discs in $\wh \C$ such that 
$B_1 \subset B_2$ and that $\partial B_1 \cap \partial B_2$ 
consists of one point $p$. 
We say that $A$ is touching at $p$. 
A closed subset in a Riemann surface homeomorphic to a crescent 
is also called a crescent. 
\end{defn}
We first consider the shape of $\LT$ in $T$: 

\begin{lem}
There exist crescents $\{{\cal A}_j\}_{j=1}^s$ and closed balls 
$\{{\cal B}_j\}_{j=1}^s$ in $T$
which satisfy the following (see Figure 7): 
\begin{enumerate} 
\item  ${\cal A}_1, \cdots , {\cal A}_s$ are mutually disjoint. 
Each ${\cal A}_j$ is homotopically equivalent to the longitude $l$ in  $T$. 
We let $p_j \in T$ denote the touching point of ${\cal A}_j$. 
\item Interiors of ${\cal B}_1, \cdots , {\cal B}_s$ 
are mutually disjoint.  
\item  The set $\bigcup_{j=1}^s {\cal B}_j$ is connected 
and homotopically equivalent to the meridian $d$ in $T$. 
In addition, it is satisfied that 
${\cal A}_j \cap \left(\bigcup_{j=1}^s {\cal B}_j \right)
=\{p_j\}$ 
for each $1 \le j \le s$, which implies that 
$\bigcup_{j=1}^s {\cal B}_j$  is a string of beads. 
\item $\LT
\subset \left(\bigcup_{j=1}^s {\cal A}_j\right) 
\cup \left(\bigcup_{j=1}^s {\cal B}_j\right)$. 
\end{enumerate}
\end{lem}
\begin{figure}[htbp]
\begin{center}
\includegraphics{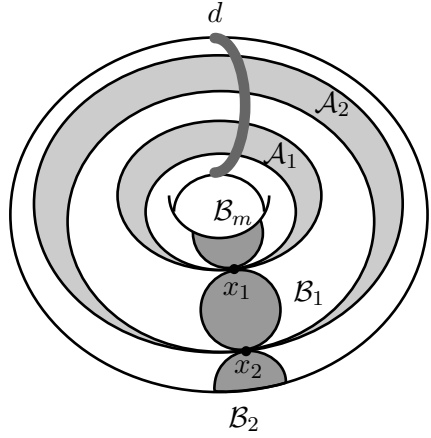}
\end{center}
\vspace{-6.5cm}
\begin{center}
\unitlength 0.1in
\begin{picture}( 40.6000, 24.5000)( 14.2000,-32.2000)
%
\special{pn 8}%
\special{pa 1420 770}%
\special{pa 5480 770}%
\special{pa 5480 3220}%
\special{pa 1420 3220}%
\special{pa 1420 770}%
\special{ip}%
\put(36.8000,-18.0000){\makebox(0,0)[lb]{${\cal A}_1$}}%
\put(39.5000,-15.3000){\makebox(0,0)[lb]{${\cal A}_2$}}%
\put(34.2000,-21.0000){\makebox(0,0)[lb]{${\cal B}_m$}}%
\put(38.4000,-25.2000){\makebox(0,0)[lb]{${\cal B}_1$}}%
\put(35.0000,-31.8000){\makebox(0,0)[lb]{${\cal B}_2$}}%
\put(33.9000,-10.4000){\makebox(0,0)[lb]{$d$}}%
\put(34.7000,-24.9000){\makebox(0,0)[lb]{$x_1$}}%
\put(35.3000,-29.0000){\makebox(0,0)[lb]{$x_2$}}%
\end{picture}%
\end{center}
\caption{Crescents $\{{\cal A}_j\}_{j=1}^s$ and 
balls $\{{\cal B}_j\}_{j=1}^s$ in $T$.}
\end{figure}
\begin{proof}
We begin with studying the subset 
$\pi_T(\Lam(\wh \Gamma) \cap \Omega_0(\Gamma_\fty))$ of $\LT$. 
Recall from Proposition 3.1 
that the set $\Lam(\wh \Gamma) \cap \Omega_0(\Gamma_\fty)$ 
is contained in the preimage 
$\pi_{Y_\infty}^{-1}(\cal N(\ul{\lam}))$
of $\rn (\ul{\lam})$ in $\Omega_0(\Gamma_\fty)$. 
Let $\wtil A$ be a connected component of 
$\pi_{Y_\infty}^{-1}(\cal N(\ul{\lam}))$ 
and let $\gamma \in \Gamma_\infty$ be a parabolic element 
such that $\gamma(\wtil A)=\wt{A}$. 
Then $\wtil A \cup \fix(\gamma)$ 
is a crescent, which mapped injectively by 
$\pi_T$ into $T$. 
If $\wtil A \cap \wtil d=\emptyset$, then 
$\wtil A \cup \fix(\gamma)$ does not separate two fixed points of $\eta$ and 
$\pi_T(\wtil A \cup \fix(\gamma))$ is contractible in $T$. 
On the other hand, 
if $\wtil A \cap \wtil d \ne \emptyset$, then 
$\wtil A \cup \fix(\gamma)$  separates two fixed points of $\eta$ and 
$\pi_T(\wtil A \cup \fix(\gamma))$ is 
homotopically equivalent to the longitude $l$ in $T$. 
Recall that $s=i(\ul{\lam}, d)$ is the intersection number 
of $\ul{\lam}$ and $d$ in $Y_\fty$.  
Then there exist exactly $s$-crescents 
$$
{\cal A}_1, \ldots, {\cal A}_s 
$$
in $T$ 
each of which is a component 
of $\pi_T \circ \pi_{Y_\infty}^{-1}(\cal N(\ul{\lam}))$ 
with its touching point $p_j \in {\cal A}_j$ 
and is homotopically equivalent to the longitude $l$.  
Since any two components of $\ul{\lam}$ are not parallel in $Y_\fty$, 
we see that $p_i \ne p_j$ for $1 \le i \ne j \le s$, and hence that 
${\cal A}_1, \cdots , {\cal A}_s$ are mutually disjoint. 
Note that the intersection  
$\LT \cap d$ lies in $\bigcup_{j=1}^s {\cal A}_j$ from 
the argument above.  
We now divide $\LT$ into 
$$
\LT=\LT^1 \sqcup \LT^2, 
$$
where 
$\LT^1=\LT \cap \bigcup_{j=1}^s {\cal A}_j$ and 
$\LT^2=\LT \setminus \bigcup_{j=1}^s {\cal A}_j$. 
Then $\LT^2$ is contained in the complement 
of $d \cup \bigcup_{j=1}^s {\cal A}_j$ in $T$, which 
consist of $s$-open balls. 
Moreover, since $\bd {\cal A}_j-\{p_j\}$ does not intersect 
$\LT$ for each $j$, 
one see that the intersection of 
the closure of $\LT^2$ with $\LT^1$ is 
$\{p_1, \ldots, p_s\}$.  
Therefore, by a slight modification of  
those open balls, we obtain a string of beads 
$\bigcup_{j=1}^s {\cal B}_j$ which contains $\LT^2$ 
and which satisfies the desired conditions. 
\end{proof}

Observe that for each $1 \le j \le s$ 
the crescent ${\cal A}_j$ contains $l_j$-crescents 
$$
{\cal A}_j^{(1)}, \ldots , {\cal A}_j^{(l_j)}
$$ 
which are components  
of $\pi_T\circ \pi_{Y_\infty}^{-1}(\cal N(\wh{\lam}))$ 
with common touching point $p_j$. 
Here $l_j$ is the weight of the component of $\lam$  
associated to ${\cal A}_j$. 
Then it is also satisfied the statement of Lemma 3.5
with ${\cal A}_j$ replaced by 
$\bigcup_{m=1}^{l_j}{\cal A}_j^{(m)}$ for all $j$. 
Adding suitable (mutually disjoint) closed neighborhoods 
$V_j$ of $p_j$,  
we obtain a submanifold 
$$
\rn_{T}=
\biggl(\, \bigcup_{j=1}^s V_j \,\biggl) \,\cup \,
\biggl(\, \bigcup_{j=1}^s \bigcup_{m=1}^{l_j}{\cal A}_j^{(m)} \,\biggl)
\,\cup\, \biggl(\, \bigcup_{j=1}^s {\cal B}_j \,\biggl)
$$
in $T$ which contains $\LT$, and a regular neighborhood 
$$
\rn(\wh \lam \cup d)=(\rn(\wh \lam) \sm d) \sqcup (\rn_{T} \sm d)   
$$
of $\wh \lam \cup d$ in $Z_\fty$ which contains 
$\LZ=(\LY \setminus d) \sqcup (\wh \Lam_{T} \setminus d)$. 
Thus we obtain the result of Proposition 3.3  
in the case of $\mu=d \in \cal S$. 
The result for general $\mu \in \mln$ is also obtained by the same argument. 

\subsection{Cutting $\rn (\wh{\lam} \cup \wh{\mu})$ 
into $\rn (\wh{(\lam,\mu)_\sharp})$ and $\rn (\wh{(\lam,\mu)_\flat})$}

Since $\Lam_{Z_n}$ converge in $Z_\fty$ to $\LZ$ 
in the sense of Lemma 2.7 
and since $\LZ$ lies in $\rn (\wh{\lam} \cup \wh{\mu})$ from Proposition 3.2,   
it follows that 
$\Lam_{Z_n}$ lie in $\rn (\wh{\lam} \cup \wh{\mu})$ 
for all large $|n|$.  
Given large enough $n>0$, 
we will show that there is a system 
of simple arcs which cuts 
$\rn (\wh{\lam} \cup \wh{\mu})$ into  
a regular neighborhood 
$\rn (\wh{(\lam,\mu)_\sharp})$ 
of a realization of $(\lam,\mu)_\sharp$ containing $\Lam_{Z_n}$; 
see Figure 8. 
The same argument yields 
$\rn (\wh{(\lam,\mu)_\flat})$ which contains $\Lam_{Z_n}$ 
for given $n<0$ with large enough $|n|$. 
\begin{figure}[htbp]
  \begin{center}
    \includegraphics{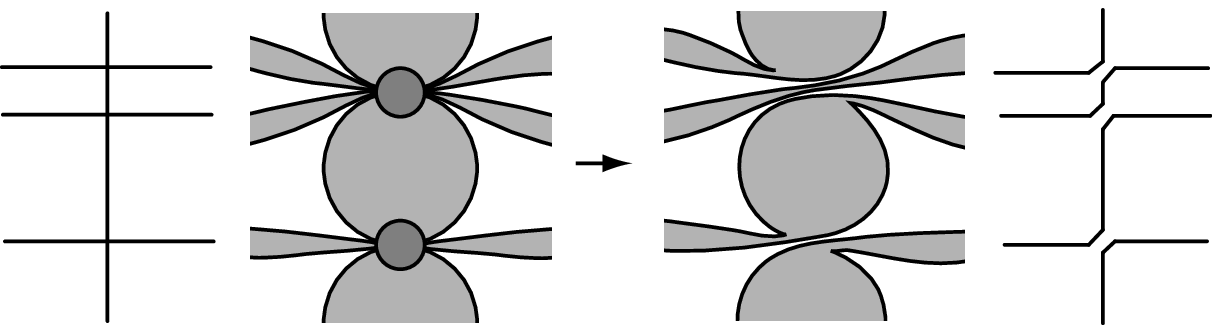}
  \end{center}
  \caption{Cutting $\rn (\wh{\lam} \cup \wh{\mu})$ 
  into $\rn (\wh{(\lam,\mu)_\sharp})$. 
The graphs in both-sides are illustrating 
$\wh{\lam} \cup \wh{\mu}$ and $\wh{(\lam,\mu)_\sharp}$ respectively. }
\end{figure}

To describe more precisely, 
we again concentrate our attention to the case of $\mu=d \in {\cal S}$. 
Let $\rn (\wh \lam \cup d)$ be a regular neighborhood 
of $\wh \lam \cup d$ containing $\LZ$ obtained as in \S 3.4. 
Note that for each $j$ the point $p_j \in \LZ$ is the pullback of 
a common fixed point $p \in \Lam(\wh \Gamma)$ 
of a maximal rank-two parabolic subgroup 
$\langle \gamma, \delta \rangle$ of $\wh \Gamma$. 
We suppose that the generators $\gamma,\,\delta$ 
satisfy the convention as in \S 3.1. 
Let $U_j \subset Z_\fty$ be an open neighborhood 
of $p_j$ containing $V_j$ and  
let $\phi_j:U_j \to U$ be a homeomorphism 
from $U_j \subset Z_\fty$ 
to an open neighborhood $U \subset \hatC$ of $p$. 
We set $V=\phi_j(V_j)$. 
Let $\omega$ be a connected component of $\Omega(\what{\Gamma})$ 
which contains $p$ in its boundary. 
Then $\omega$ is simply connected and is 
invariant under the action of $\gamma$. 
We may assume that $\omega \cap U$ 
consists of exactly two components $\omega_1$ and $\omega_2$ 
for which 
$\gamma(\omega_1) \subset \omega_1$ 
and $\gamma^{-1}(\omega_2) \subset \omega_2$; see Figure 9. 
In addition, we assume that   
$\delta(\omega_2)$ and $\delta^{-1}(\omega_2)$ 
are contained in $U$. 
In this situation, we have the following: 

\begin{prop}
Let $x \in \omega_1, \, y \in \omega_2$ 
such that $x,\,y,\,\delta(y), \delta^{-1}(y) \not\in V$. 
Then 
there exist simple arcs $\beta_n$ 
 which 
satisfy the following: 
\begin{enumerate}
\item $\beta_n \subset \Omega(\Gamma_n)$ for all $|n| \gg 0$, 
\item $\beta_n$ joins $x$ and $\delta(y)$ and 
is contained in 
$\omega_1 \cup V \cup \delta(\omega_2)$ for all $n \gg0$, and   
\item $\beta_n$ joins $x$ and $\delta^{-1}(y)$ and 
is contained in 
$\omega_1 \cup V \cup \delta^{-1}(\omega_2)$ for all $n \ll 0$.  
\end{enumerate}
\end{prop}
\begin{figure}[htbp]
\begin{center}
\includegraphics{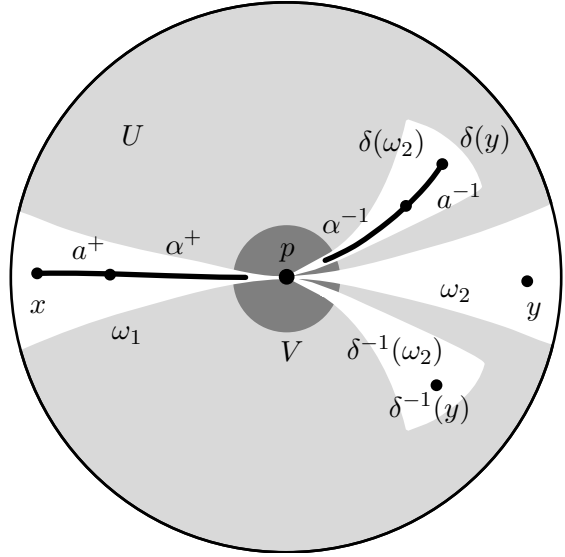}
\end{center}
\vspace{-8.5cm}
\begin{center}
\unitlength 0.1in
\begin{picture}( 37.6000, 30.9000)( 11.6000,-35.9000)
%
\special{pn 8}%
\special{pa 1160 500}%
\special{pa 4920 500}%
\special{pa 4920 3590}%
\special{pa 1160 3590}%
\special{pa 1160 500}%
\special{ip}%
\put(21.9000,-13.6000){\makebox(0,0)[lb]{$U$}}%
\put(17.1000,-22.4000){\makebox(0,0)[lb]{$x$}}%
\put(19.3000,-19.6000){\makebox(0,0)[lb]{$a^+$}}%
\put(24.2000,-19.3000){\makebox(0,0)[lb]{$\alpha^+$}}%
\put(21.4000,-24.0000){\makebox(0,0)[lb]{$\omega_1$}}%
\put(38.6000,-21.8000){\makebox(0,0)[lb]{$\omega_2$}}%
\put(43.1000,-22.8000){\makebox(0,0)[lb]{$y$}}%
\put(33.7000,-25.2000){\makebox(0,0)[lb]{$\delta^{-1}(\omega_2)$}}%
\put(35.9000,-28.2000){\makebox(0,0)[lb]{$\delta^{-1}(y)$}}%
\put(39.6000,-14.1000){\makebox(0,0)[lb]{$\delta(y)$}}%
\put(38.4000,-16.9000){\makebox(0,0)[lb]{$a^{-1}$}}%
\put(32.4000,-18.2000){\makebox(0,0)[lb]{$\alpha^{-1}$}}%
\put(30.2000,-19.7000){\makebox(0,0)[lb]{$p$}}%
\put(30.2000,-24.9000){\makebox(0,0)[lb]{$V$}}%
\put(34.3000,-14.3000){\makebox(0,0)[lb]{$\delta(\omega_2)$}}%
\end{picture}%
\end{center}
\caption{Figure explaining Proposition 3.7 and its proof. }
\end{figure}

Before proving the proposition above, 
We shall complete the proof of Theorem 1.3. 
We fix $1 \le j \le s$ for a wile. 
Recall that $p_j \in U_j$ is the touching point 
of $l_j$-crescents $A_j^{(1)}, \ldots , A_j^{(l_j)}$ in $T$.  
Thus we may assume that $U_j \setminus \rn(\wh \lam \cup d)$ 
consists of exactly $2l_j+2$ connected components. 
We first consider the case of $n \gg 0$. 
Given large enough $n>0$,  
it follows from Proposition 3.5 that 
there are mutually disjoint $l_j$-simple arcs 
$$
\beta_n^{(1)}, \ldots , \beta_n^{(l_j)}
$$
in $U_j$ which satisfy the following: 
\begin{itemize}
\item $\beta_n^{(m)}$ does not intersect $\LZ$ for each $1 \le m \le  l_j$, 
\item the intersection 
$\beta_n^{(m)} \cap \rn(\wh \lam \cup d)$ is contained in $V_j$ 
for each $1 \le m \le l_j$, and 
\item all $2l_j$-end points of $\beta_n^{(1)}, \ldots , \beta_n^{(l_j)}$ 
are contained in distinct components of 
$U_j \setminus \rn(\wh \lam \cup d)$. 
\end{itemize}
Observe that 
by cutting $\rn(\wh \lam \cup d)$ 
along those arcs contained in $U_j$ for every $1 \le j \le s$, 
we obtain a regular neighborhood 
$\rn(\wh{(\lam,d)_\sharp})$ of a realization of $(\lam,d)_\sharp$ 
which contains $\Lam_{Z_n}$; 
see Figure 8. 
Note that since $Z_n$ lies in $Q(S)$, 
Lemma 2.6 implies that $\Lam_{Z_n}$ consists of mutually 
disjoint non-trivial simple closed curves 
contained in $\rn(\wh{(\lam,d)_\sharp})$. 
We claim that 
each connected component of $\rn(\wh{(\lam,d)_\sharp})$
contains exactly two connected components of $\Lam_{Z_n}$, 
and hence that $\Lam_{Z_n}$ is a realization of $2 (\lam, d)_{\sharp}$. 
In fact, the claim holds because 
each connected component of $\rn(\wh{(\lam,d)_\sharp})$
contains a connected component of 
$\rn(\wh \lam) \setminus d$ which contains 
exactly two essential connected components of $\LY \setminus d$; 
see \cite{It1} for more rigorous argument. 
Thus Theorem 2.6 implies that $Z_n \in \Q_{(\lam,d)_\sharp}$. 
Similarly, we have $Z_n  \in  \Q_{(\lam,d)_\flat}$ for all $n \ll 0$. 
Therefore we have completed the proof of Theorem 1.3 for the case of 
$\mu=d \in {\cal S}$. 
The result for general $\mu \in \mln$ is obtained by the same argument. 
All what we have to do is to show Proposition 3.5. 

\begin{proof}[Proof of Proposition 3.5]
We only consider the case of $n \gg 0$, 
since the argument for the case $n \ll 0$ is completely parallel.  
Recall that there is a sequence $\gamma_n \in \Gamma_n$ such that 
$\gamma_n \to \gamma$ and $(\gamma_n)^n \to \delta$ as $n \to +\infty$. 
Assume for simplicity that $y$ lies in the 
$\langle \gamma \rangle$-orbit of $x$; 
i.e. $y=\gamma^{-n_0}(x)$ for some positive integer $n_0$. 
Then $\delta(y)=\gamma^{-n_0} \delta(x)$ holds since $\gamma$ commutes with 
$\delta$. 
Let $N>0$ be a positive integer such that 
for all $k \ge N$, 
$\gamma^k(x)$ and $\gamma^{-k} \delta(y)$ are contained in $V$. 
Let $a^+$ be an arc in $\omega_1$ joining $x$ and $\gamma(x)$. 
Since $\Omega(\Gamma_n)$ converge to $\Omega(\what{\Gamma})$ 
in the sense of Carath\'{e}odory (see \cite{KeTh}),  
$a^+ \subset \Omega(\wh{\Gamma})$ are contained in $\Omega(\Gamma_n)$ 
for all $n \gg 0$.   
Hence there exist arcs $b_n^+$ joining 
$x$ and $\gamma_n(x)$ in $\Omega(\Gamma_n)$ 
which satisfy $b_n^+ \hto a^+$ as $n \to +\infty$. 
Now set 
\begin{eqnarray*}
\beta_n=\bigcup_{k=0}^{n-n_0-1}(\gamma_n)^k(b_n^+). 
\end{eqnarray*}
Then $\beta_n$ is an arc which joins $x$ to $(\gamma_n)^{n-n_0}(x)$. 
Since $(\gamma_n)^{n-n_0}(x)$ converge to 
$\gamma^{-n_0}\delta(x)=\delta(y)$ as $n \to +\fty$, 
we may assume that the arc $\beta_n$ joins $x$ to $\delta(y)$. 
We will show that $\beta_n$ satisfy the desired conditions. 
Since $b_n^+ \subset \Omega(\Gamma_n)$ and $\gamma_n \in \Gamma_n$, 
it follows that $\beta_n \subset \Omega(\Gamma_n)$. 
Let us now consider two subarcs 
\begin{eqnarray*}
\beta_n^+=\bigcup_{k=0}^{N}(\gamma_n)^k(b_n^+), \quad 
\beta_n^-=\bigcup_{k=n-n_0-1-N}^{n-n_0-1}(\gamma_n)^k(b_n^+) 
\end{eqnarray*}
of $\beta_n$, which are the unions of the first $N$-orbits and 
the last $N$-orbits, respectively. 
Then we will obtain $\beta_n \subset \omega_1 \cup V \cup \delta(\omega_2)$ 
by showing that $\beta_n^+ \subset \omega_1$, 
$\beta_n^- \subset \delta(\omega_2)$ 
and $\beta_n-(\beta_n^+ \cup \beta_n^-) \subset V$ 
for all $n \gg 0$. 

\bigskip
{\bf Claim 1}: \ 
 $\beta_n^+ \subset \omega_1$ and $\beta_n^- \subset \delta(\omega_2)$ 
 for all $n \gg 0$. 

\bigskip

Let $a^-$ denote the arc $\gamma^{-n_0-1}\delta(a^+)$ 
joining $\gamma^{-1}\delta(y)$ and $\delta(y)$. We may 
assume that $a^- \subset \delta(\omega_2)$.  
Let us consider two arcs 
\begin{eqnarray*}
\alpha^+=\bigcup_{k=0}^{N}\gamma^k(a^+), \quad 
\alpha^-=\bigcup_{k=0}^{N}\gamma^{-k}(a^-), 
\end{eqnarray*}
which are contained in $\omega_1$ and $\delta(\omega_2)$, respectively. 
To obtain the claim, 
it suffices to show that 
$\beta_n^+ \hto \alpha^+$ and $\beta_n^- \hto \alpha^-$ in $\wh \C$. 
Since $\gamma_n \to \gamma$ and $b_n^+ \hto a^+$, 
it follows that $\beta_n^+ \hto \alpha^+$. 
Now set $b_n^-=(\gamma_n)^{n-n_0-1}(b_n^+)$. 
Then we have $b_n^- \hto a^-=\gamma^{-n_0-1}\delta(a^+)$ 
from $(\gamma_n)^{n-n_0-1} \to \gamma^{-n_0-1}\delta$ 
and $b_n^+ \hto a^+$.  
Thus we conclude that 
\begin{eqnarray*}
\beta_n^-=\bigcup_{l=0}^{N}(\gamma_n)^{-l}(\gamma_n)^{n-n_0-1}(b_n^+)
=\bigcup_{l=0}^{N}(\gamma_n)^{-l}(b_n^-)
 \ \hto \  
\bigcup_{l=0}^{N}(\gamma_n)^{-l}(a^-)=\alpha^-. 
\end{eqnarray*}

\bigskip
{\bf Claim 2}: \ 
$\beta_n-(\beta_n^+ \cup \beta_n^-) \subset V$ for all $n \gg 0$. 
\bigskip

We will show that the orbits 
$\{(\gamma_n)^k(x)\,|\,N<k<n-n_0-1-N\}$ of $x$ 
in $\beta_n-(\beta_n^+ \cup \beta_n^-)$ 
are contained in $V$. 
The claim is obtained 
by the same argument. 
Since $V$ is a neighborhood of the common fixed point $p$ of 
$\langle  \gamma, \delta \rangle$, 
we may assume that 
by changing $N$ large enough if necessary, 
$|s| \le N$ holds whenever 
$\gamma^s \delta^t(x)$ is not contained in the interior $\inte(V)$ of $V$. 
To obtain a contradiction, suppose that 
there is a sequence $k_n \to \infty \,(n \to \infty)$ 
such that $N<k_n<n-n_0-1-N$ and that 
$(\gamma_n)^{k_n}(x) \not\in V$. 
Then passing to a subsequence if necessary, 
$(\gamma_n)^{k_n}(x)$ converges to some point 
$\hat{x} \not\in \inte(V)$. 
Then the sequence $(\gamma_n)^{k_n}$ does not diverge in $\psl$ 
because both fixed points of $\gamma_n$ converge to $p$ and 
$(\gamma_n)^{k_n}(x)$ converge to $\hat{x} \ne p$. 
Since 
$\langle \gamma_n \rangle \hto \langle \gamma, \delta \rangle$ in $\psl$, 
after passing to a further subsequence, 
we may assume that 
$(\gamma_n)^{k_n} \to \gamma^s \delta^t$ 
for some $s,\,t \in \Z$. 
Then it follows that 
$\hat{x}=\gamma^s \delta^t(x) \not\in \inte(V)$, and hence that 
$|s| \le N$. 
On the other hand, since the convergence 
$(\gamma_n)^{k_n} \to \gamma^s \delta^t$ implies the convergence 
$(\gamma_n)^{k_n-tn-s} \to \id$ as $n \to +\fty$, 
we have 
$k_n \equiv tn+s \, (n \gg 0)$ from Lemma 3.6 in \cite{JoMa}. 
It then follows 
from $N<k_n<n-n_0-1-N$ that $t=0, \,s \ge 0$ or $t=1,\,s \le 0$, 
both of which 
contradict to $|s| \le N$. 
\end{proof}

\section{Corollaries of the main theorem}

We collect here some consequences of Theorem 1.3. 

\begin{thm}
For any non-zero $\lambda \in \mln$,
there exists 
$Y \in \ov{\Q_0} \cap \ov{\Q_\lam}$ 
such that $U \cap \Q_{\lambda}$ is disconnected 
for any sufficiently small neighborhood $U$ of $Y$.  
In particular, $\ov{\Q_\lam}$ is not a topological manifold with
boundary.
\end{thm}

\begin{proof}
The result follows directly from Theorem 1.3 
in the case of $i(\lam, \mu) \ne 0$. 
\end{proof}

We will show in Theorem 4.3 that any two components of $Q(S)$ bump 
from Theorem 1.3 combined with 
the following result, which is a generalization of Theorem 1.2.

\begin{thm}[Theorem B in \cite{It1}]
Let $\{\lam_i\}_{i=1}^m$ be a finite subset of $\mln-\{0\}$ such that 
$i(\lam_i,\lam_j)=0$ for any $1 \le i < j \le m$. 
Then we have 
$\ov{\Q_0} \cap \ov{\Q_{\lam_1}} \cap \cdots \cap 
\ov{\Q_{\lam_m}} \ne \emptyset$. 
\end{thm}

\begin{thm}
For any $\lam, \mu \in \mln$, we have 
$\overline{\Q_\lam} \cap \overline{\Q_\mu} \ne\emptyset$.
\end{thm}

\begin{proof}
Let $\lam,\,\mu \in \mln$.  
If $i(\lam,\mu)=0$, then the result follows from Theorem 4.2. 
From now on, 
we assume that $i(\lam,\mu) \ne 0$ 
and decompose $\mu$ into 
$\mu=\mu'+\mu''$ so that $\mu', \, \mu'' \in \mln$ and 
$i(\lambda,\mu)=i(\lambda,\mu')$ are satisfied. 
The proof is divided into two cases 
(i) $\mu''=0$ and (ii) $\mu'' \ne 0$. 

(i) 
Case of $\mu''=0$. Then $\mu=\mu'$. 
By the same argument as in \S 3, replacing $\lam$ to $(\lam,\mu)_\flat$, 
we obtain a sequence 
$$
\{Y_n\}_{|n| \gg 0} \subset \Q_{(\lambda,\mu)_{\flat}}
$$ 
which converges to some $Y_\fty \in \bd \Q_0$  
as $|n| \to \infty$.  
Since $(\lambda,\mu)_{\flat}$ and $\mu$ have no parallel component 
in common, we obtain a grafting 
$\Gr_{\mu}(Y_{\infty}) \in \bd \Q_\mu$ of $Y_\fty$ 
and a convergent sequence 
$$
Z_n \to \Gr_\mu(Y_\fty) \quad (|n| \to \infty)
$$
which satisfies $\hol(Z_n)=\hol(Y_n)$ for all large $|n|$. 
It follows from Theorem 1.3 that 
$Z_n \in \Q_{((\lambda,\mu)_{\flat},\mu)_{\sharp}}=\Q_{\lambda}$ 
for all $n \gg 0$. 
Thus we obtain 
$\overline{\Q_\lam} \cap \overline{\Q_\mu} \ne\emptyset$.

(ii) 
Case of $\mu'' \ne 0$. 
Since $i(\lam,\mu'')=0$ and $i(\mu', \mu'')=0$,  
we have $i((\lam, \mu')_{\flat}, \mu'')=0$. 
Theorem 4.2 then implies that 
$
\ov{\Q_0} \cap \ov{\Q_{(\lam, \mu')_\flat}} \cap \ov{\Q_{\mu''}} \ne \emptyset.
$
More precisely, we obtain as in \S 3 two sequences
$$
\{Y_n\}_{|n| \gg 0} \subset \Q_{(\lambda, \mu')_{\flat}}, 
\quad \{Y'_n\}_{|n| \gg 0} \subset  \Q_{\mu''}
$$ 
both of which converge to some $Y_{\infty} \in \partial \Q_0$ 
as $|n| \to \infty$; see the proof of Theorem B in \cite{It1} 
for more details. 
Note that the holonomy image of $Y_\fty$ is a $b$-group 
whose parabolic locus is the support of $(\lambda,\mu')_{\flat}+\mu''$. 
Since $(\lambda,\mu')_{\flat}+\mu''$ and $\mu'$ 
have no parallel component in common, 
we obtain a grafting $\Gr_{\mu'}(Y_{\infty})$ 
and two convergent sequences 
$$
Z_n, \,Z_n'\to \Gr_{\mu'}(Y_\fty) \quad (|n| \to \infty)
$$
which satisfy $\hol(Z_n)=\hol(Y_n)$ and 
$\hol(Z'_n)=\hol(Y'_n)$ for all large $|n|$. 
It follows from Theorem 1.3 that 
$Z_n \in \Q_{(( \lam, \mu' )_\flat, \mu')_\sharp}=\Q_\lam$
and $Z'_n \in \Q_{\mu' + \mu''} =\Q_\mu$
for all $n \gg 0$. 
Thus we obtain 
$\ov{\Q_\lam} \cap \ov{\Q_\mu} \ne\emptyset$ also in this case. 
\end{proof}

\begin{thm}
For any non-zero $\lam \in \mln$, 
the holonomy map $\hol:P(S) \to R(S)$ 
is not injective on $\ov{\Q_\lam}$, 
although it is invective on $\Q_\lam$. 
\end{thm}

\begin{proof} 
Let $\nu$ be an element of $\mln$ each of whose 
component intersects $\lam$ essentially.  
Lemma 2.3 then implies that 
\begin{eqnarray*}
(\nu,(\lam,\nu)_\sharp)_\sharp=((\lam,\nu)_\sharp,\nu)_\flat=\lam, \\
(\nu,(\lam,\nu)_\flat)_\flat=((\lam,\nu)_\flat,\nu)_\sharp=\lam. 
\end{eqnarray*} 
Now let  
$\{Y_n\}_{|n| \gg 0} \subset Q_\nu$ be a sequence 
constructed as in \S 3 which converges to some 
$Y_\fty \in \bd \Q_0$ as $|n| \to \fty$. 
Applying Theorem 1.3, we obtain two sequences in $\Q_\lam$ as follows: 
\begin{itemize}
\item $\{Z_n\}_{n \gg 0}$ in 
$\Q_{(\nu,(\lam,\nu)_\sharp)_\sharp}=\Q_\lam$ which converges to 
$\Gr_{(\lam,\nu)_\sharp}(Y_\fty)$ as $n \to +\fty$. 
\item $\{Z'_n\}_{n \ll 0}$ in 
$\Q_{(\nu,(\lam,\nu)_\flat)_\flat}=\Q_\lam$ which converges to 
$\Gr_{(\lam,\nu)_\flat}(Y_\fty)$ as $n \to -\fty$. 
\end{itemize}
Since $i(\lam,\nu) \ne 0$, we have 
$(\lam,\nu)_\sharp \ne (\lam,\nu)_\flat$. 
Therefore  
$\Gr_{(\lam,\nu)_\sharp}(Y_\fty) \ne \Gr_{(\lam,\nu)_\flat}(Y_\fty)$ 
and the result follows. 
\end{proof}

\bigskip
\begin{flushleft}
Graduate School of Mathematics, \\
Nagoya University, \\
Nagoya 464-8602, Japan  \\
\texttt{itoken@math.nagoya-u.ac.jp} 
\end{flushleft}


\begin{thebibliography}{ABCD} 

\bibitem[AC]{AnCa}J.~W.~Anderson and R.~D.~Canary, 
  {\it Algebraic limits of Kleinian groups which rearrange the pages of a book}, 
  Invent. Math. {\bf 126} (1996), 205--214. 

\bibitem[ACM]{AnCaMc}J.~W.~Anderson, R.~D.~Canary and D.~Metallurgy, 
{\it On the topology of deformation spaces of Kleinian groups}, 
Ann. of Math. {\bf 152} (2000), 693--741. 

\bibitem[Be1]{Be1}L.~Bers, 
{\it Simultaneous uniformization}, 
Bull. Amer. Math. Soc. {\bf 66} (1960), 94--97. 

\bibitem[Be2]{Be2}L.~Bers, 
{\it On boundaries of Theichm\"{u}ller spaces and on Kleinian groups, I}, 
Ann. of math. {\bf 91} (1970), 570--600. 


\bibitem[BO]{BoOt}F.~Bonahon and J.~P.~Otal, 
{\it Vari\'{e}t\'{e}s hyperboliques \`{a} 
g\'{e}od\'{e}siques arbitrairement courtes}, 
Bull. London Math. Soc. {\bf 20} (1988), 255--261. 

\bibitem[BH]{BrHo}K.~Bromberg and J.~Holt, 
{\it Self-bumping of deformation spaces of hyperbolic 3-manifolds}, 
J. Diff. Geom. {\bf 57} (2001), 47--65. 

\bibitem[Ca]{Ca}R.~D.~Canary, 
{\it Pushing the boundary}, 
in the Tradition of Ahlfors and Bers, III, Contemp. Math. 
{\bf 355}, Amer. Math. Soc., 2004, 109--121. 

\bibitem[Co]{Co}T.~D.~Comar, 
{\it Hyperbolic Dehn surgery and convergence of Kleinian groups}, 
Ph. D. Thesis, University of Michigan, 1996. 

\bibitem[Ea]{Ea}C.~J.~Earle, 
{\it On variation of projective structures}, 
Riemann surfaces and related topics, 
Ann. Math. Studies {\bf 97} (1981), 87--99. 

\bibitem[Ga]{Ga}D.~M.~Gallo, 
{\it Deforming real projective structures}, 
Ann. Acad. Sci. Fenn. {\bf 22} (1997), 3--14. 

\bibitem[Go]{Go}W.~M.~Goldman, 
  {\it Projective structures with Fuchsian holonomy}, 
  J. Diff.  Geom. {\bf 25} (1987), 297--326. 

\bibitem[He]{He}D.~A.~Hejhal, 
{\it Monodromy groups and linearly polymorphic functions}, 
Acta Math.  {\bf 135} (1975), 1--55. 

\bibitem[Ho]{Ho}J.~Holt, 
{\it Bumping and self-bumping of deformation spaces}, 
in the Tradition of Ahlfors and Bers, III, Contemp. Math. 
{\bf 355}, Amer. Math. Soc., 2004, 269--284.  

\bibitem[Hu]{Hu}J.~H.~Hubbard, 
{\it The monodromy of projective structures}, 
Riemann surfaces and related topics, 
Ann. Math. Studies {\bf 97} (1981), 257--275. 

\bibitem[It1]{It1}K.~Ito, 
{\it Exotic projective structures and quasi-Fuchsian space}, 
Duke Math. J. {\bf 105} (2000), 185--209. 

\bibitem[It2]{It2} K. ~Ito, 
{Grafting and components of quasi-fuchsian projective structures}, 
to appear in Spaces of Kleinian groups   
(eds. Y.~Minsky, M.~Sakuma, C.~Series), London Math. Soc. 
Lecture note series {\bf 329}, Cambridge University Press. 

\bibitem[Ita]{Ita}K.~Ito, 
{\it On continuous extensions of grafting maps}, 
preprint. 
\texttt{math.GT/0411133}

\bibitem[JM]{JoMa}T.~J{\o}rgensen and A.~Marden, 
{\it Algebraic and geometric convergence of Kleinian groups}, 
Math. Scand. {\bf 66} (1990), 47--72. 

\bibitem[KS]{KoSu}Y.~Komori and T.~Sugawa, 
{\it Bers embedding of the Teichmuller space of a once-punctured torus}, 
Conform. Geom. Dyn. {\bf 8} (2004), 115--142.  

\bibitem[KSYW]{KoSuYaWa}Y.~Komori, T.~Sugawa, Y.~Yamashita, M.~Wada,  
{\it Drawing Bers embeddings of the Teichmuller space of once-punctured tori}, 
to appear in Experimental Mathematics. 

\bibitem[KT]{KeTh}S.~P.~Kerckhoff and W.~P.~Thurston, 
  {\it Non-continuity of the action of the modular group 
   at Bers' boundary of Teichmuller space}, 
  Invent. Math. {\bf 100} (1990), 25--47. 
  
\bibitem[Lu]{Lu}F.~Luo, 
{\it Some applications of a multiplicative structure on simple loops in surfaces}, in 
 Knots, braids, and mapping class groups---papers dedicated to Joan S. Birman (New York, 1998), 
AMS/IP Stud. Adv. Math., {\bf 24}, Amer. Math. Soc., Providence, 2001, 123--129. 

\bibitem[Ma]{Ma}A.~Marden, 
{\it The geometry of finitely generated kleinian groups}, 
Ann. of Math. {\bf 99} (1974), 383--462. 

\bibitem[MT]{MaTa}K.~Matsuzaki and M.~Taniguchi, 
 {\it Hyperbolic manifolds and Kleinian groups}, 
Oxford University Press, 1998. 

\bibitem[Mc]{Mc}C.~T.~McMullen, 
  {\it Complex earthquakes and Teichmuller theory}, J. Amer. Math. Soc. 
{\bf 11} (1998), 283--320. 

\bibitem[Su]{Su}D.~P.~Sullivan, 
{\it Quasiconformal homeomorphisms and dynamics II: 
structural stability implies hyperbolicity of Kleinian groups}, 
Acta Math. {\bf 155} (1985), 243--260. 

\end{thebibliography}
\end{document}